\documentclass{amsart}
\usepackage{amsmath,amssymb}
\usepackage{yhmath}

\newcommand{\bdis}{\begin{displaymath}}
\newcommand{\edis}{\end{displaymath}}
\newcommand{\be}{\begin{equation}}
\newcommand{\ee}{\end{equation}}
\newcommand{\mbb}{\mathbb}
\newcommand{\mcal}{\mathcal}

\newcommand{\vp}{\varphi}

\newcommand{\zf}{\zeta\left(\frac{1}{2}+it\right)}

\theoremstyle{definition}

\theoremstyle{remark}
\newtheorem{remark}[]{Remark}

\newtheorem*{mydef11}{{\bf Theorem 1}}

\newtheorem*{mydef12}{{\bf Theorem 2}}

\newtheorem*{mydef13}{{\bf Theorem 3}}

\newtheorem*{mydef14}{{\bf Theorem 4}}

\newtheorem*{mydef15}{{\bf Theorem 5}}

\newtheorem*{mydef21}{{\bf Definition 1}}

\newtheorem*{mydef22}{{\bf Definition 2}}

\newtheorem*{mydef23}{{\bf Definition 3}}

\newtheorem*{mydef41}{{\bf Corollary 1}}

\newtheorem*{mydef42}{{\bf Corollary 2}}

\newtheorem*{mydef51}{{\bf Lemma 1}}

\newtheorem*{mydef52}{{\bf Lemma 2}}

\newtheorem*{mydef53}{{\bf Lemma 3}}

\newtheorem*{mydef54}{{\bf Lemma 4}}

\newtheorem*{mydef55}{{\bf Lemma 5}}

\newtheorem*{mydef56}{{\bf Lemma 6}}

\newtheorem*{mydef81}{{\bf Property 1}}

\newtheorem*{mydef82}{{\bf Property 2}}

\numberwithin{equation}{section}

\begin{document}

\title{Jacob's ladders, almost linear increments of the Hardy-Littlewood integral (1918) and their relations to the Selberg's formula (1946) and the Fermat-Wiles theorem}{}

\author{Jan Moser}

\address{Department of Mathematical Analysis and Numerical Mathematics, Comenius University, Mlynska Dolina M105, 842 48 Bratislava, SLOVAKIA}

\email{jan.mozer@fmph.uniba.sk}


\begin{abstract}
In this paper we give new consequences that follow from our formula for increments of the Hardy-Littlewood integral. Main of these consequences is an $\zeta$-equivalent of the Fermat-Wiles theorem. It is expressed purely by means of the Riemann's zeta-function on the critical line and by the Jacob's ladder.  
\end{abstract}
\maketitle

\section{Introduction} 

\subsection{}    

Let us remind that in our paper \cite{7} we have obtained new result of completely new type for certain class of increments of the classical Hardy-Littlewood integral\footnote{See \cite{1}.}
\be \label{1.1} 
J(T)=\int_0^T\left|\zf\right|^2{\rm d}t. 
\ee 
Namely: 
\be \label{1.2} 
\begin{split}
& \int_{\overset{r-1}{T}}^{\overset{r}{T}}\left|\zf\right|^2{\rm d}t=(1-c)\overset{r-1}{T}+\mcal{O}(T^{a+\delta}),\\ 
& a\in \left[\frac 14,\frac 13\right],\ r=1,2,\dots,k 
\end{split}
\ee 
for every fixed $k\in\mbb{N}$, where 
\be \label{1.3} 
\overset{r}{T}=\vp_1^{-r}(T),\ \forall\- T>T_0>0,\ \overset{0}{T}=T, 
\ee 
and $c$ is the Euler's constant and $T_0$ is sufficiently big positive number. 

Next, the symbol $a$ denotes the least valid exponent for which the estimate of the error term in the Hardy-Littlewood-Ingham formula holds true, see \cite{7}, (2.5), (2.7), $0<\delta$ being sufficiently small. 

\begin{remark}
For the beginning it is sufficient to put $a=\frac 13$. 
\end{remark} 

\begin{remark}
The almost linear formula in variable $\overset{r-1}{T}$ represents newly discovered property of the Hardy-Littlewood integral (\ref{1.1}) on the continuum set of segments 
\be \label{1.4} 
\begin{split}
& [\overset{r-1}{T},\overset{r}{T}]:\ \overset{r}{T}-\overset{r-1}{T}\sim (1-c)\pi(\overset{r}{T}), \\ 
& \overset{r}{T}=\overset{r}{T}(T),\ T\to\infty,\ r=1,\dots,k, 
\end{split}
\ee 
(see \cite{7}, (3.3)), where 
\be \label{1.5} 
\pi(\overset{r}{T})\sim \frac{\overset{r}{T}}{\ln\overset{r}{T}},\ T\to\infty 
\ee 
and the symbol $\pi(x)$ stands for the prime-counting function. 
\end{remark} 

\subsection{} 

In this paper we use the following notions: 
\begin{itemize}
	\item[{\tt (a)}] Jacob's ladder $\vp_1(t)$, 
	\item[{\tt (b)}] the function 
	\bdis 
	\begin{split}
	& \tilde{Z}^2(t)=\frac{{\rm d}\vp_1(t)}{{\rm d}t}=\frac{1}{\omega(t)}\left|\zf\right|^2,\\ 
	& \omega(t)=\left\{1+\mcal{O}\left(\frac{\ln\ln t}{\ln t}\right)\right\}\ln t,\ t\to\infty, 
	\end{split}
	\edis 
	\item[{\tt (c)}] direct iterations of Jacob's ladders 
	\bdis 
	\begin{split}
	& \vp_1^0(t)=t,\ \vp_1^1(t)=\vp_1(t),\ \vp_1^2(t)=\vp_1(\vp_1(t)),\dots , \\ 
	& \vp_1^k(t)=\vp_1(\vp_1^{k-1}(t))
	\end{split}
	\edis 
	for every fixed natural number $k$, 
	\item[{\tt (d)}] reverse iterations of Jacob's ladders 
	\bdis 
	\begin{split}
	& \vp_1^{-1}(T)=\overset{1}{T},\ \vp_1^{-2}(T)=\vp_1^{-1}(\overset{1}{T})=\overset{2}{T},\dots, \\ 
	& \vp_1^{-r}(T)=\vp_1^{-1}(\overset{r-1}{T})=\overset{r}{T},\ r=1,\dots,k, 
	\end{split} 
	\edis 
	where 
	\bdis 
	\begin{split}
	& \overset{0}{T}=T<\overset{1}{T}(T)<\overset{2}{T}(T)<\dots<\overset{k}{T}(T), \\ 
	& T\sim \overset{1}{T}\sim \overset{2}{T}\sim \dots\sim \overset{k}{T},\ T\to\infty. 
	\end{split}
	\edis  
\end{itemize}  

We have introduced all these notions into the theory of the Riemann's zeta-function in our works \cite{2} -- \cite{5}. 

\subsection{} 

Next, let us remind the following two functions: 
\be \label{1.6} 
S(t)=\frac{1}{\pi}\arg\zf, 
\ee 
where the argument is defined by the usual way\footnote{Comp. \cite{8}, p. 90.} and 
\be \label{1.7} 
S_1(t)=\int_0^t S(\tau){\rm d}\tau. 
\ee 
We have found new results concerning functions $S(t)$ and $S_1(t)$. The first type result is represented by the formula 
\be \label{1.8} 
\begin{split}
& \int_T^{\overset{1}{T}}\left|\zf\right|^2{\rm d}t+\frac{1-c}{d(l)}\int_T^{\overset{k-1}{T}}|S_1(t)|^{2l}{\rm d}t\sim \\ 
&  \int_{\overset{k-1}{T}}^{\overset{k}{T}}\left|\zf\right|^2{\rm d}t,\ T\to\infty, 
\end{split}
\ee 
where $d(l)$ is the Selberg's constant, see \cite{8}, p. 130:\ $c_k\to d(l)$ for every fixed $l\in\mbb{N}$. The second type of result is then represented by the formula: 
\be \label{1.9} 
\begin{split}
& \frac{1}{\overset{r}{T}}\int_{\overset{r-1}{T}}^{\overset{r}{T}}\left\{ d(l)\left|\zf\right|^2+(1-c)|S_1(t)|^{2l}\right\}{\rm d}t\sim \\ 
& (1-c)d(l),\ r=1,\dots,k,\ T\to\infty. 
\end{split}
\ee 

\begin{remark}
\emph{The lift transformation} $\mcal{T}_{0,k-1}$: 
\be \label{1.10} 
\int_T^{\overset{1}{T}}\left|\zf\right|^2{\rm d}t\xrightarrow{\mcal{T}_{0,k-1}}\int_{\overset{k-1}{T}}^{\overset{k}{T}}\left|\zf\right|^2{\rm d}t 
\ee 
is defined by the formula (\ref{1.8}) by means of the integral 
\be \label{1.11} 
\frac{1-c}{d(l)}\int_T^{\overset{k-1}{T}}|S_1(t)|^{2l}{\rm d}t 
\ee 
as a needful portion of \emph{energy} for realization of \emph{transmutation} (\ref{1.10}). 
\end{remark} 

\begin{remark}
Our formula (\ref{1.9}) can be viewed as a kind of \emph{the conservation law}\footnote{Comp. \cite{7}, subsection 4.1.} for very complex oscillatory dynamics generated by the following linear combinations 
\be \label{1.12} 
\begin{split}
& d(l)\left|\zf\right|^2+(1-c)|S_1(t)|^{2l},\ t\in [\overset{r-1}{T},\overset{r}{T}],\ r=1,\dots,k,\ T\to\infty. 
\end{split}
\ee 
\end{remark} 

\subsection{} 

Let us remind that the equation 
\be \label{1.13} 
x^n+y^n=z^n
\ee 
is by the Fermat - Wiles theorem not solvable in positive integers for any $n\geq 3$. 

In this paper we obtain new result for the set 
\be \label{1.14}
\left\{\frac{x^n+y^n}{z^n}\right\},\ \forall\- x,y,z,n\in\mbb{N}
\ee 
of Fermat's rationals. 

\begin{remark}
The $\zeta$-condition 
\be \label{1.15} 
\lim_{\tau\to 0}\frac{1}{\tau}\int_{\frac{x^n+y^n}{z^n}\frac{\tau}{1-c}}^{\left[\frac{x^n+y^n}{z^n}\frac{\tau}{1-c}\right]^1}\left|\zf\right|^2{\rm d}t\not=1 
\ee 
on the set (\ref{1.14}) represents the $\zeta$-equivalent of the Fermat-Wiles theorem. 
\end{remark} 

\begin{remark}
Consequently, the $\zeta$-condition (\ref{1.15}) can be viewed as the first point contact between the almost linear increments (\ref{1.2}) of the Hardy-Littlewood integral and the Fermat-Wiles theorem. More generally, between the Riemann's zeta-function on the critical line and the Fermat-Wiles theorem. 
\end{remark} 

\section{On the point of contact between Riemann's zeta-function on the critical line and the set of $L_2$-orthogonal systems} 

In this section we remind our preceeding results, see \cite{6}, since these give new point of contact with the theory of $L_2$ orthogonal systems that is completely different field of analysis. 

\subsection{} 

We have introduced generating vector operator $\hat{G}$ defined on the class of all $L_2$-orthogonal systems 
\be \label{2.1} 
\{f_n(t)\}_{n=0}^\infty,\ t\in [a,a+2l],\ \forall \- a\in\mbb{R},\ \forall \- l>0
\ee 
by the following action: 
\be \label{2.2} 
\begin{split}
& \{f_n(t)\}_{n=0}^\infty\xrightarrow{\hat{G}}\{f^{p_1}_n(t)\}_{n=0}^\infty\xrightarrow{\hat{G}}\{f^{p_1,p_2}_n(t)\}_{n=0}^\infty\xrightarrow{\hat{G}}\dots \\ 
& \xrightarrow{\hat{G}}\{f^{p_1,p_2,\dots,p_s}_n(t)\}_{n=0}^\infty,\ p_1,\dots,p_s=1,\dots,k 
\end{split}
\ee  
for every fixed $k,s\in\mbb{N}$ with the explicit formula\footnote{See \cite{6}, (2.19)} for 
\bdis 
f^{p_1,p_2,\dots,p_s}_n(t). 
\edis 

\subsection{} 

For example, in the case of Legendre $L_2$-orthogonal system 
\be \label{2.3} 
\{ P_n(t)\}_{n=0}^\infty,\ t\in [-1,1] 
\ee 
the operator $\hat{G}$ produces the third generation as follows: 
\be \label{2.4} 
\begin{split}
& P_n^{p_1,p_2,p_3}(t)=P_n(u_{p_1}(u_{p_2}(u_{p_3}(t))))\times\prod_{r=0}^{p_1-1}\left|\tilde{Z}(v_{p_1}^r(u_{p_2}(u_{p_3}(t))))\right| \times \\ 
& \prod_{r=0}^{p_2-1}\left|\tilde{Z}(v_{p_2}^r(u_{p_3}(t)))\right|\times \prod_{r=0}^{p_3-1}\left|\tilde{Z}(v_{p_3}(t))\right|, \\ 
& p_1,p_2,p_3=1,\dots,k,\ t\in[-1,1],\ a=-1,\ l=1, 
\end{split} 
\ee  
where 
\be \label{2.5} 
u_{p_i}(t)=\vp_1^{p_i}\left(\frac{\overset{p_i}{\wideparen{T+2}}-\overset{p_i}{T}}{2}(t+1)+\overset{p_i}{T}\right)-T-1,\ i=1,2,3
\ee 
are automorphisms on $[-1,1]$ and 
\be \label{2.6} 
\begin{split}
& v_{p_i}^r(t)=\vp_1^{p_i}\left(\frac{\overset{p_i}{\wideparen{T+2}}-\overset{p_i}{T}}{2}(t+1)+\overset{p_i}{T}\right),\ r=0,1,\dots,p_i-1, \\ 
& t\in [-1,1] \ \Rightarrow \ u_{p_1}(t)\in [-1,1] \ \wedge \ v_{p_i}^r(t)\in [\overset{p_i-r}{T},\overset{p_i-r}{\wideparen{T+2}}]. 
\end{split}
\ee 

\begin{mydef81}
\begin{itemize}
	\item[(a)] Every member of every new $L_2$-orthogonal system 
	\bdis 
	\{ P_n^{p_1,p_2,p_3(t)}\}_{n=0}^\infty,\ t\in [-1,1],\ p_1,p_2,p_3=1,\dots,k 
	\edis  
	contains the function 
	\bdis 
	\left|\zf\right|_{t=\tau}
	\edis 
	for correpsonding $\tau$ since\footnote{See \cite{3}, (9.1), (9.2).} 
	\be \label{2.7} 
	|\tilde{Z}(t)|=\sqrt{\frac{{\rm d}\vp_1(t)}{{\rm d}t}}=\{ 1+o(1)\}\frac{1}{\sqrt{\ln t}}\left|\zf\right|,\ t\to\infty. 
	\ee 
	\item[(b)] Property (a) holds true due to Theorem of the paper \cite{6} for every generation 
	\bdis 
	\{ f_n^{p_1,\dots,p_s}(t)\}_{n=0}^\infty,\ t\in [a,a+2l], s\in\mbb{N}. 
	\edis 
\end{itemize}
\end{mydef81} 

\begin{remark}
We see that there is a close binding between the theory of the Riemann's zeta-function on the critical line and the theory of $L_2$-orthogonal systems. 
\end{remark} 

\section{The point of contact between the increments of the Hardy-Littlewood integral and the Selberg's formulas}

\subsection{}  

Let us remind: 

\begin{itemize}
	\item[{\tt(a)}] Selberg's formulas \cite{8}, p. 130: 
	\be \label{3.1} 
	\int_{T}^{T+H}|S_1(t)|^{2l}{\rm d}t=d(l)H+\mcal{O}\left(\frac{H}{\ln T}\right), 
	\ee 
	\be \label{3.2} 
	T^b\leq H\leq T,\ \frac 12<b\leq 1,\ d(l)>0, 
	\ee 
	\be \label{3.3} 
	\int_0^T|S_1(t)|^{2l}{\rm d}t=d(l)T+\mcal{O}\left(\frac{T}{\ln T}\right). 
	\ee 
	\item[{\tt(b)}] Our almost linear formula (\ref{1.2}) 
	\be \label{3.4} 
	\int_{\overset{r-1}{T}}^{\overset{r}{T}}\left|\zf\right|^2{\rm d}t=(1-c)\overset{r-1}{T}+\mcal{O}(T^{a+\delta}),\ r=1,\dots,l 
	\ee 
	in the context of formulae (\ref{3.1}) -- (\ref{3.3}). 
\end{itemize} 

\subsection{} 

By making use of the substitution 
\be \label{3.5} 
T\to \overset{r-1}{T},\ H\to \overset{r}{T}-\overset{r-1}{T},\ r=1,\dots,k 
\ee 
in (\ref{3.1}), we obtain 
\be \label{3.6} 
\int_{\overset{r-1}{T}}^{\overset{r}{T}}|S_1(t)|^{2l}{\rm d}t=d(l)(\overset{r}{T}-\overset{r-1}{T})+\mcal{O}\left(\frac{\overset{r}{T}-\overset{r-1}{T}}{\overset{r-1}{T}}\right). 
\ee  
Next, by the formulae (\ref{1.4}) and (\ref{1.5}) and by the subsection 1.2 
\be \label{3.7} 
\begin{split}
& T\sim\overset{1}{T}\sim \dots\sim \overset{k}{T}, \\ 
& \overset{r}{T}-\overset{r-1}{T}\sim(1-c)\frac{\overset{r}{T}}{\ln \overset{r}{T}}\sim (1-c)\frac{T}{\ln T}
\end{split}
\ee 
and it follows that 
\be \label{3.8} 
H=\overset{r}{T}-\overset{r-1}{T},\ \mcal{O}\left(\frac{H}{\ln T}\right)=\mcal{O}\left(\frac{T}{\ln^2 T}\right). 
\ee 
Consequently we obtain the following statement. 

\begin{mydef51}
\be \label{3.9} 
\frac{1}{d(l)}\int_{\overset{r-1}{T}}^{\overset{r}{T}}|S_1(t)|^{2l}{\rm d}t=\overset{r}{T}-\overset{r-1}{T}+\mcal{O}\left(\frac{T}{\ln^2 T}\right) 
\ee  
for every fixed $k,l\in\mbb{N}$. 
\end{mydef51} 

Further, the summation of the formula 
\be \label{3.10} 
\frac{1}{d(l)}\int_{\overset{n-1}{T}}^{\overset{n}{T}}|S_1(t)|^{2l}{\rm d}t=\overset{n}{T}-\overset{n-1}{T}+\mcal{O}\left(\frac{T}{\ln^2 T}\right) 
\ee  
over the segment $[r,s-1]\subset [1,k]$ implies the following: 
\be \label{3.11} 
\frac{1}{d(l)}\int_{\overset{r-1}{T}}^{\overset{s-1}{T}}|S_1(t)|^{2l}{\rm d}t=\overset{s-1}{T}-\overset{r-1}{T}+\mcal{O}\left(\frac{T}{\ln^2 T}\right). 
\ee 

Now, our formulas (\ref{3.4}) and (\ref{3.11}) give us the Theorem 1. 

\begin{mydef11}
For every fixed $k,l\in\mbb{N}$ it is true that 
\be \label{3.12} 
\begin{split}
& \int_{\overset{s-1}{T}}^{\overset{s}{T}}\left|\zf\right|^2{\rm d}t= \\ 
& \int_{\overset{r-1}{T}}^{\overset{r}{T}}\left|\zf\right|^2{\rm d}t+\frac{1-c}{d(l)}\int_{\overset{r-1}{T}}^{\overset{s-1}{T}}|S_1(t)|^{2l}{\rm d}t+\mcal{O}\left(\frac{T}{\ln^2 T}\right), 
\end{split}
\ee  
where 
\be \label{3.13} 
1\leq r\leq s-1\leq k-1. 
\ee 
\end{mydef11} 

\begin{remark}
\emph{The lift transformation} 
\be \label{3.14} 
\int_{\overset{r-1}{T}}^{\overset{r}{T}}\left|\zf\right|^2{\rm d}t\xrightarrow{\mcal{T}_{r-1,s-1}}\int_{\overset{s-1}{T}}^{\overset{s}{T}}\left|\zf\right|^2{\rm d}t
\ee 
is defined by (\ref{3.12}). This lift transformation acts on the set 
\be \label{3.15} 
\left\{
\int_{T}^{\overset{1}{T}}\left|\zf\right|^2{\rm d}t,\ \int_{\overset{1}{T}}^{\overset{2}{T}}\left|\zf\right|^2{\rm d}t,\dots,\int_{\overset{k-1}{T}}^{\overset{k}{T}}\left|\zf\right|^2{\rm d}t
\right\}
\ee 
of integrals by means of the integral 
\be \label{3.16} 
\frac{1-c}{d(l)}\int_{\overset{r-1}{T}}^{\overset{s-1}{T}}|S_1(t)|^{2l}{\rm d}t, 
\ee  
from the level $r-1$ onto the level $s-1$. 
\end{remark} 

\begin{remark}
The integral (\ref{3.16}) can be viewed as needful portion of energy that is absorbed by the integral 
\bdis 
\int_{\overset{r-1}{T}}^{\overset{r}{T}}\left|\zf\right|^2{\rm d}t
\edis  
to realize \emph{the transmutation} (\ref{3.14}). 
\end{remark} 

\subsection{} 

There is also another possibility. Namely, the sum of formula (\ref{3.9}) and the following one\footnote{See (\ref{3.4}).} 
\be \label{3.17} 
\frac{1}{1-c}\int_{\overset{r-1}{T}}^{\overset{r}{T}}\left|\zf\right|^2{\rm d}t=\overset{r-1}{T}+\mcal{O}(T^{a+\delta})
\ee 
give us the formula 
\be \label{3.18} 
\begin{split}
& \int_{\overset{r-1}{T}}^{\overset{r}{T}}\left\{
\frac{1}{1-c}\left|\zf\right|^2+\frac{1}{d(l)}|S_1(t)|^{2l}
\right\}{\rm d}t=\overset{r}{T}+\mcal{O}\left(\frac{T}{\ln^2 T}\right). 
\end{split}
\ee  

Consequently we have the following Theorem. 

\begin{mydef12}
\be \label{3.19} 
\begin{split}
& \frac{1}{\overset{r}{T}}\int_{\overset{r-1}{T}}^{\overset{r}{T}}
\left\{ 
d(l)\left|\zf\right|^2+(1-c)|S_1(t)|^{2l}
\right\}{\rm d}t=\\ 
& (1-c)d(l)+\mcal{O}\left(\frac{1}{\ln^2T}\right),\ T\to\infty,\ r=1,\dots,k. 
\end{split}
\ee 
\end{mydef12} 

\subsection{} 

Next, if we use the substitution 
\be \label{3.20} 
T\to \overset{r-1}{T},\ r=1,\dots,k 
\ee  
in our formula (\ref{3.3}) then we obtain\footnote{See (\ref{3.8}).} 
\be \label{3.21} 
\int_0^{\overset{r-1}{T}}|S_1(t)|^{2l}{\rm d}t=d(l)\overset{r-1}{T}+\mcal{O}\left(\frac{T}{\ln T}\right). 
\ee  
Now, we get Theorem 3 by (\ref{3.17}) and (\ref{3.21}). 

\begin{mydef13}
It is true that 
\be \label{3.22} 
\int_0^{\overset{r-1}{T}}|S_1(t)|^{2l}{\rm d}t=\frac{d(l)}{1-c}\int_{\overset{r-1}{T}}^{\overset{r}{T}}\left|\zf\right|^2{\rm d}t+\mcal{O}\left(\frac{T}{\ln T}\right), 
\ee  
and, of course, 
\be \label{3.23} 
\int_{\overset{r-1}{T}}^{\overset{r}{T}}\left|\zf\right|^2{\rm d}t=\frac{1-c}{d(l)}\int_0^{\overset{r-1}{T}}|S_1(t)|^{2l}{\rm d}t+\mcal{O}\left(\frac{T}{\ln T}\right), 
\ee 
where 
\be \label{3.24} 
[0,\overset{r-1}{T})\cap (\overset{r-1}{T},\overset{r}{T}]=\emptyset,\ r=1,\dots,k. 
\ee 
\end{mydef13} 

\begin{remark}
The following effect is expressed by the formulae (\ref{3.22}) and (\ref{3.23}): each element of the 2-tuple 
\be \label{3.25} 
\left\{ 
\int_{\overset{r-1}{T}}^{\overset{r}{T}}\left|\zf\right|^2{\rm d}t,\ \int_0^{\overset{r-1}{T}}|S_1(t)|^{2l}{\rm d}t
\right\} ,\ r=1,\dots,k 
\ee  
of integrals generates asymptotically the complementary element. This effect can be compared with nonlinear formulae we have obtained in \cite{7}, subsection 6.1. 
\end{remark} 

\section{New $\zeta$-functional generated by almost linear formula (\ref{1.2})} 

\subsection{} 

In the following we shall use our fomula (\ref{1.2}) with $r=1$, that is 
\be \label{4.1} 
\int_T^{\overset{1}{T}}\left|\zf\right|^2{\rm d}t=(1-c)T+\mcal{O}(T^{a+\delta}),\ T>T_0>0, 
\ee  
where $T_0$ is sufficiently big and\footnote{See subsection 1.2, {\tt(d)}.} 
\be \label{4.2} 
\overset{1}{T}=[T]^1=\vp_1^{-1}(T). 
\ee 
Now, if we put 
\be \label{4.3} 
T=\frac{x}{1-c}\tau,\ \tau\in \left(\frac{1-c}{x}T_0,+\infty\right),\ x>0 
\ee 
into (\ref{4.1}), then we obtain the following statement. 

\begin{mydef52}
\be \label{4.4} 
\frac{1}{\tau}\int_{\frac{x}{1-c}\tau}^{\left[\frac{x}{1-c}\tau\right]^1}\left|\zf\right|^2{\rm d}t=x+\mcal{O}(\tau^{-1+a+\delta}), 
\ee  
where 
\be \label{4.5} 
\tau\in (\tau_1(x),+\infty),\ \tau_1(x)=\max\left\{ \left(\frac{1-c}{x}\right)^2,(T_0)^2\right\},\ x>0 
\ee  
and the constant in the $\mcal{O}$-term depends on $x$. 
\end{mydef52} 

Consequently, we have the following. 

\begin{mydef53}
\be \label{4.6} 
\lim_{\tau\to\infty}\frac{1}{\tau}\int_{\frac{x}{1-c}\tau}^{\left[\frac{x}{1-c}\tau\right]^1}\left|\zf\right|^2{\rm d}t=x 
\ee  
for every fixed $x>0$, that is the formula (\ref{4.6}) associates to every fixed ray 
\be \label{4.7} 
y(\tau;x)=\frac{x}{1-c}\tau,\ \tau\in(\tau_1(x),+\infty) 
\ee 
with given positive $x$, that is defined by the slope 
\be \label{4.8} 
0<\arctan\frac{x}{1-c}<\frac{\pi}{2}. 
\ee 
\end{mydef53} 

\begin{mydef21}
We shall call the mapping $\mcal{F}_1$: 
\be \label{4.9} 
y(\tau;x)\xrightarrow{\mcal{F}_1}x 
\ee 
defined by (\ref{4.6}) as $\zeta$-functional. 
\end{mydef21} 

\begin{remark}
Let us remind that the elements $y(\tau;x)$ are defined on different neighbourhoods of plus infinity, see (\ref{4.7}). 
\end{remark} 

\subsection{} 

It is clear that the formula (\ref{4.6}) implies the following. 

\begin{mydef82}
If $x_1,\dots,x_n>0$, then: 
\be \label{4.10} 
\begin{split}
& \lim_{\tau\to\infty}\frac{1}{\tau}\int_{\frac{\sum x_i}{1-c}\tau}^{\left[\frac{ \sum x_i}{1-c}\tau\right]^1}\left|\zf\right|^2{\rm d}t= \\ 
& \sum_{i=1}^n \lim_{\tau\to\infty}\frac{1}{\tau}\int_{\frac{x_i}{1-c}\tau}^{\left[\frac{x_i}{1-c}\tau\right]^1}\left|\zf\right|^2{\rm d}t, 
\end{split}
\ee  
\be \label{4.11} 
\begin{split}
& \lim_{\tau\to\infty}\frac{1}{\tau}\int_{\frac{\prod x_i}{1-c}\tau}^{\left[\frac{ \prod x_i}{1-c}\tau\right]^1}\left|\zf\right|^2{\rm d}t= \\ 
& \prod_{i=1}^n \lim_{\tau\to\infty}\frac{1}{\tau}\int_{\frac{x_i}{1-c}\tau}^{\left[\frac{x_i}{1-c}\tau\right]^1}\left|\zf\right|^2{\rm d}t, 
\end{split}
\ee  
\be \label{4.12} 
\begin{split}
& \lim_{\tau\to\infty}\frac{1}{\tau}\int_{\frac{x_1/x_2}{1-c}\tau}^{\left[\frac{x_1/x_2}{1-c}\tau\right]^1}\left|\zf\right|^2{\rm d}t= \\ 
& \lim_{\tau\to\infty}
\frac
{\int_{\frac{x_1}{1-c}\tau}^{\left[\frac{x_1}{1-c}\tau\right]^1}\left|\zf\right|^2{\rm d}t} 
{\int_{\frac{x_2}{1-c}\tau}^{\left[\frac{x_2}{1-c}\tau\right]^1}\left|\zf\right|^2{\rm d}t}, 
\end{split}
\ee 
where factors $\frac{1}{\tau}$ cancel in the last expression. 
\end{mydef82} 

\section{On the existence of the point of contact between the almost linear increments (\ref{1.2}) of the Hardy-Littlewood integral and the Fermat-Wiles theorem}

If we use the substitution 
\be \label{5.1} 
x\to \frac{x^n+y^n}{z^n}, 
\ee  
the we obtain the following. 

\begin{mydef54}
\be \label{5.2} 
\lim_{\tau\to\infty}\frac{1}{\tau}\int_{\frac{x^n+y^n}{z^n}\frac{\tau}{1-c}}^{\left[\frac{x^n+y^n}{z^n}\frac{\tau}{1-c}\right]^1}\left|\zf\right|^2{\rm d}t=\frac{x^n+y^n}{z^n}  
\ee 
for every fixed Fermat's rational. 
\end{mydef54} 

Consequently, we have the following. 

\begin{mydef14}
The $\zeta$-condition 
\be \label{5.3} 
\lim_{\tau\to\infty}\frac{1}{\tau}\int_{\frac{x^n+y^n}{z^n}\frac{\tau}{1-c}}^{\left[\frac{x^n+y^n}{z^n}\frac{\tau}{1-c}\right]^1}\left|\zf\right|^2{\rm d}t\not= 1
\ee  
on the class of all Fermat's rationals represents the $\zeta$-equivalent of the Fermat-Wiles theorem. 
\end{mydef14} 

Next, if we use the formula (\ref{4.12}) in the case (\ref{5.3}) the we obtain the following. 

\begin{mydef41} 
The $\zeta$-condition 
\be \label{5.4} 
\lim_{\tau\to\infty}
\frac
{\int_{(x^n+y^n)\frac{\tau}{1-c}}^{\left[(x^n+y^n)\frac{\tau}{1-c}\right]^1}\left|\zf\right|^2{\rm d}t}
{\int_{z^n\frac{\tau}{1-c}}^{\left[z^n\frac{\tau}{1-c}\right]^1}\left|\zf\right|^2{\rm d}t}\not= 1 
\ee  
represents the second $\zeta$-equivalent of the Fermat-Wiles theorem. 
\end{mydef41}  

\section{Next kinds of $\zeta$-functionals and $\zeta$-equivalents of the Fermat-Wiles theorem} 

\subsection{} 

It is true by (\ref{4.1}) and Remark 1 that 
\be \label{6.1} 
\begin{split}
& \int_T^{\overset{1}{T}}\left|\zf\right|^2{\rm d}t=(1-c)T+\mcal{O}(T^{1/3+\delta})= \\ 
& T(1-c)\{ 1+\mcal{O}(T^{-2/3+\delta})\}= \\ 
& T\exp\left[ \ln((1-c)\{1+\mcal{O}(T^{-2/3+\delta})\})\right]. 
\end{split}
\ee 
Putting 
\be \label{6.2} 
T=x^\tau 
\ee 
into (\ref{6.1}) we obtain: 

\begin{mydef55}
\be \label{6.3} 
\lim_{\tau\to\infty}\left\{
\int_{x^\tau}^{[x^\tau]^1}\left|\zf\right|^2{\rm d}t
\right\}^{\frac{1}{\tau}}=x, 
\ee 	 
where 
\be \label{6.4} 
x>1,\ \tau\in (\tau_2(x),+\infty),\ \tau_2(x)=\max\left\{\frac{1}{\ln^2x},\ln^2T_0\right\}, 
\ee 
i.e. the formula (\ref{6.3}) associates the value $x$ to every fixed exponential function $x^\tau$. 
\end{mydef55} 

\begin{mydef22}
	We shall call the mapping 
	\be \label{6.5} 
	x^\tau\xrightarrow{\mcal{F}_2}x,\ x>1 
	\ee  
	as the second $\zeta$-functional. 
\end{mydef22} 

\subsection{} 

Next, if we put 
\be \label{6.6} 
T=\tau^x,\ x>0,\ \tau\in((T_0)^{1/x},+\infty) 
\ee 
into (\ref{6.1}), we obtain, that 
\be \label{6.7} 
\begin{split}
& \frac{1}{\ln\tau}\ln\left\{\int_{\tau^x}^{[\tau^x]^1}\left|\zf\right|^2{\rm d}t\right\}= \\ 
& x+\frac{1}{\ln\tau}\ln[(1-c)\{1+o(1)\}]. 
\end{split}
\ee 
As a simple consequence we have the following lemma. 

\begin{mydef56}
\be \label{6.8} 
\lim_{\tau\to\infty}\ln\left\{\int_{\exp[x\ln\tau]}^{[\exp[x\ln\tau]]^1}\left|\zf\right|^2{\rm d}t\right\}^{\frac{1}{\ln\tau}}=x, 
\ee  
where 
\be \label{6.9} 
\tau\in(\tau_3(x),+\infty),\ \tau_3(x)=\max\left\{(T_0)^{1/x},T_0\right\},\ x>0, 
\ee  
i.e. the formula (\ref{6.8}) associates the value $x$ to every fixed power function $\tau^x$. 
\end{mydef56} 

\begin{mydef23}
	We shall call the mapping $\mcal{F}_3$: 
	\be \label{6.10} 
	\tau^x\xrightarrow{\mcal{F}_3}x,\ x>0 
	\ee  
	as the third $\zeta$-functional.\footnote{Enumeration of $\zeta$-functionals is related with this paper only.}
\end{mydef23} 

\subsection{} 

If we use the substitution 
\be \label{6.11} 
x\to \frac{x^n+y^n}{z^n}
\ee 
in (\ref{6.8}), then we obtain the following statement (comp. section 5). 

\begin{mydef15}
The $\zeta$-condition 
\be \label{6.12} 
\lim_{\tau\to\infty}\ln\left\{
\int_{\exp(\frac{x^n+y^n}{z^n}\ln \tau)}^{[\exp(\frac{x^n+y^n}{z^n}\ln \tau)]^1}\left|\zf\right|^2{\rm d}t
\right\}^{\frac{1}{\ln\tau}}\not=1 
\ee  
on the class of all Fermat's rationals represents the third $\zeta$-equivalent of the Fermat-Wiles theorem. 
\end{mydef15} 

\begin{remark}
Of course, there are analogues of the properties (\ref{4.10}) -- (\ref{4.12}) that are true for functionals (\ref{6.3}) and (\ref{6.8}). 
\end{remark} 

Consequently, the $\zeta$-condition (\ref{6.12}) implies the following. 

\begin{mydef42}
\be \label{6.13} 
\lim_{\tau\to\infty}
\frac
{\ln\{\int_{\exp((x^n+y^n)\ln\tau)}^{[\exp((x^n+y^n)\ln\tau)]^1}|\zf|^2{\rm d}t\}}
{\ln\{\int_{\exp(z^n\ln\tau)}^{[\exp(z^n\ln\tau)]^1}|\zf|^2{\rm d}t\}}
\not=1
\ee 
represents the fourth $\zeta$-equivalent of the Fermat-Wiles theorem, where 
\bdis 
[\exp((x^n+y^n)\ln\tau)]^1=\vp_1^{-1}(\tau^{x^n+y^n}),\dots 
\edis 
and $\vp_1^{-1}$ is the first reverse iteration of the Jacob's ladder. 
\end{mydef42} 

I would like to thank Michal Demetrian for his moral support of my study of Jacob's ladders.


\begin{thebibliography}{29}
\bibitem{1} 
G.H. Hardy, J.E. Littlewood, Contribution to the theory of the Riemann zeta-function and the theory of the distribution of Primes, Acta Math. 41 (1), 119 -- 196, (1918). 
\bibitem{2}
J. Moser,
`Jacob's ladders and almost exact asymptotic representation of the Hardy-Littlewood integral`,
Math. Notes 88, (2010), 414-422, arXiv: 0901.3937. 
\bibitem{3}
J. Moser,
`Jacob's ladders, the structure of the Hardy-Littlewood integral and some new class of nonlinear integral equations`,
Proc. Steklov Inst. 276 (2011), 208-221, arXiv: 1103.0359. 
\bibitem{4}
J. Moser, Jacob's ladders, reverse iterations and new infinite set of $L_2$-orthogonal systems generated by the Riemann $\zf$-function, arXiv: 1402.2098v1. 
\bibitem{5} 
J. Moser, Jacob's ladders, interactions between $\zeta$-oscillating systems and $\zeta$-analogue of an elementary trigonometric identity, Proc. Stek. Inst. 299, 189 -- 204, (2017). 
\bibitem{6} 
J. Moser, Jacob's ladders and vector operator producing new generations of $L_2$-orthogonal systems connected with the Riemann's $\zf$ function, arXiv: 2302.0750v3. 
\bibitem{7} 
J. Moser, Jacob's ladders, existence of almost linear increments of the Hardy-Littlewood integral and new types of multiplicative laws, arXiv: 2304.09267. 
\bibitem{8} 
A. Selberg. Contributions to the theory of the Riemann zeta function, Arch. f\" ur Math. og Naturv. B, 48, pp. 89 -- 155, (1946). 
\end{thebibliography}
\end{document}